\documentclass[12pt]{amsart}

\usepackage{amssymb}
\usepackage{epsf}

\advance \topmargin by -\headheight
\advance \topmargin by -\headsep
     
\evensidemargin \oddsidemargin
\advance\hoffset by -3mm  
\advance\voffset by  8mm  

\theoremstyle{plain}
\newtheorem{theorem}{Theorem}

\theoremstyle{remark}
\newtheorem{remark}{Remark}

\newcommand{\oone}{\hbox{${\mathcal O}_1$}}
\newcommand{\otwo}{\hbox{${\mathcal O}_2$}}
\newcommand{\othr}{\hbox{${\mathcal O}_3$}}
\newcommand{\othrd}{\hbox{${\mathcal O}_3'$}}
\newcommand{\ofou}{\hbox{${\mathcal O}_4$}}

\begin{document}

\title{Kirby Calculus in Manifolds with Boundary}
\author{Justin Roberts}
\address{Department of Mathematics\\University of California\\
Berkeley, CA 94720\\USA}
\curraddr{Department of Mathematics and Statistics\\University of
Edinburgh\\Edinburgh EH9 3JZ\\Scotland}
\email{justin@maths.ed.ac.uk}
\date{January 18, 1996}
\maketitle

Suppose there are two framed links in a compact, connected 3-manifold
(possibly with boundary, or non-orientable) such that the associated
3-manifolds obtained by surgery are homeomorphic (relative to their
common boundary, if there is one.) How are the links related? Kirby's
theorem \cite{K1} gives the answer when the manifold is $S^3$, and
Fenn and Rourke \cite{FR} extended it to the case of any closed
orientable 3-manifold, or $S^1 \tilde{\times} S^2$. The purpose of
this note is to give the answer in the general case, using only minor
modifications of Kirby's original proof.

Let $M$ be a compact, connected, orientable (for the moment)
3-manifold with boundary, containing (in its interior) a framed link
$L$. Doing surgery on this link produces a new manifold, whose
boundary is canonically identified with the original $\partial M$. In
fact {\em any} compact connected orientable $N$, whose boundary is
identified (via some chosen homeomorphism) with that of $M$, may be
obtained by surgery on $M$ in such a way that the boundary
identification obtained after doing the surgery agrees with the chosen
one. This is because $M \cup (\partial M \times I) \cup N$ (gluing $N$
on via the prescribed homeomorphism of boundaries) is a closed
orientable 3-manifold, hence bounds a (smooth orientable) 4-manifold,
by Lickorish's theorem \cite{L1}. Taking a handle decomposition of
this 4-manifold starting from a collar $M \times I$ requires no
0-handles (by connectedness) and no $1$- or $3$-handles, because
these may be traded (surgered 4-dimensionally) to 2-handles (see
\cite{K1}). The attaching maps of the remaining 2-handles determine a
framed link $L$ in $M$, surgery on which produces $N$.

\begin{figure}
\[ \vcenter{\epsffile{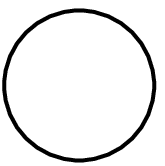}}\  \longrightarrow\ \vcenter{\epsffile{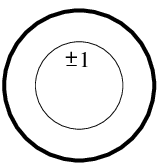}}\qquad \vcenter{\epsffile{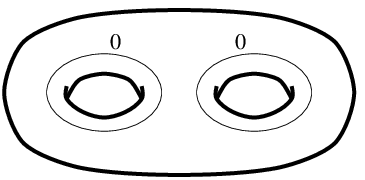}}\
\longrightarrow\ \vcenter{\epsffile{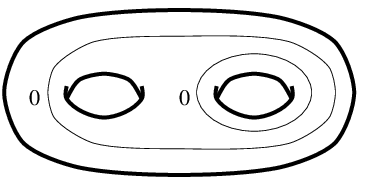}}\]
\caption{{\em Blow-up} (\oone) and {\em Handleslide} (\otwo) moves. \label{f:12}}
\end{figure}

The framed link representation is not at all unique, and the natural
question is: given framed links $L_0$ and $L_1$ in $M$ such that the
surgered manifolds $M_0, M_1$ are homeomorphic relative to their
boundary (there is a canonical identification between these boundaries
which we must not change), how are $L_0$ and $L_1$ related? If $M$ is
the 3-sphere, the answer was given by Kirby \cite{K1}: there is a
finite sequence of (isotopy classes of) links, the first being $L_0$
and the last $L_1$, such that each is obtained from its predecessor by
a move of type
\oone\ or \otwo\ or its reverse. The move \oone\ is supported in a
3-ball in $M$: it is simply disjoint union with a $\pm 1$-framed
unknot.  The move \otwo\ is supported in a genus-2 handlebody in
$M$: it is any embedded image of the pattern depicted in figure
\ref{f:12}, which is a modification of zero-framed
links occurring inside a standard unknotted handlebody in $S^3$.
(This is probably easier than thinking about it as a
parallel-and-connect-sum operation.)

\begin{figure}
\[ \vcenter{\epsffile{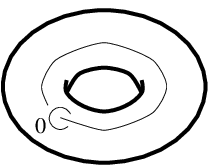}}\ \longrightarrow\
\vcenter{\epsffile{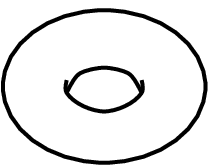}}
\]
\caption{{\em Circumcision} (\othr) move.\label{f:3}}
\end{figure}

Fenn and Rourke proved in \cite{FR} that the theorem could be extended
to any closed orientable $M$ by allowing an additional move \othr,
supported in a solid torus, and shown in figure \ref{f:3}; its
amusing name is due to Kauffman. (It is well-known that in $S^3$
this move follows from \oone\ and \otwo.)

\begin{theorem}
For an arbitrary compact, connected, orientable 3-manifold with boundary, moves
\oone, \otwo, \othr\ still suffice.
\end{theorem}

\begin{proof}
The proof follows from a careful reading of the original proof in
Kirby's paper \cite{K1}: little additional effort is required.
However, to convince the reader, a short summary of this proof,
recalling the main stages and noting the differences which arise, will
be presented.

{\em Stage 1.} Let $W_0$ and $W_1$ be 4-manifolds formed by taking $M
\times I$ and attaching 2-handles to the top surface along the links
$L_0$, $L_1$. They have common boundaries $\partial W_0 = \partial W_1
= M \cup (\partial M \times I) \cup N$, so cross this with $I$ and use
it to glue them together to get a closed smooth 4-manifold. After
connect-summing $W_0$ with some `${\mathbb C}P^2$'s or their reverses
(move \oone) the signature of this closed manifold may be taken to be
zero, so it bounds a smooth 5-manifold $V$. Form a handle
decomposition of $V$, founded on $W_0 \times I$, by taking a generic
Morse function on $V$ which restricts to 0 on $W_0$, 1 on $W_1$, and
$t$ on the slice $\partial W_0 \times \{t\}$ of the connecting product
4-manifold. This decomposition can be chosen to have no 0- or
5-handles (by connectedness and non-closedness), and all its handles'
attaching maps land in $W_0$, because of the form of the Morse
function on the boundary. Just as in \cite{K1}, surger the 1- and
4-handles to 3- and 2-handles, reorder them in the natural way and
view $V$ as having two ``ends'' $W_0$ and $W_1$, and a ``middle
slice'' $W_{1/2}$ (the level set of $\frac12$) which is reached,
working inwards from a collar on either end, by the attachment of
5-dimensional 2-handles. Such things attach by framed `$S^1 \times
B^3$'s lying in whichever end $W_i$ ($i=0,1$) is being
considered. Since each $W_i$ is obtained by attaching just 2-handles
to $M \times I$, the images of the attaching maps may be isotoped into
$M \times I$. Now comes a crucial difference: in \cite{K1}, each $S^1
\times B^3$ lies in $S^3
\times I$, so it is null-homotopic and may be unknotted; the effect of
surgery is thus to connect-sum $W_i$ with $(S^4-S^1
\times B^3) \cup (B^2 \times S^2)$, which is either $S^2 \times S^2$
or the twisted bundle $S^2
\tilde{\times} S^2$, depending on which of the two framings of $S^1
\times B^3$ is used. Such a connect-sum can be achieved by using only
\oone\ and \otwo\ (although of course it may be regarded as just the 
reverse of \othr). In the general situation, the homotopy class of the
knot in $M \times I$ may be non-trivial, and it is clear that the
attaching cannot in general be achieved using just \oone\ and \otwo,
because these moves preserve the subgroup of $\pi_1(M)$ generated by
the components of the framed link. Instead, if the $S^1 \times B^3$ is
pictured as lying on top of the top surface of $M \times I$, meeting
it in $S^1 \times B^2$ (some framed knot whose ``integral'' framing
stabilised to a ``mod 2'' one for $S^1 \times B^3$ agrees with the
given one), then replacing it by $B^2 \times S^2$ may be achieved by
decomposing $B^2 \times S^2 = B^2
\times B^2 \cup_{B^2 \times S^1} B^2 \times B^2$ and adding these
pieces one at a time: the first attaches as a 2-handle to the framed
knot $S^1 \times B^2$, and the second as a 2-handle to the framed core
of that handle. Pushing it off the core into the top surface of $M$,
results in a small zero-framed meridional curve linking the original
one once. This is exactly the reverse of the move \othr, occurring in
the framed solid torus $S^1 \times B^2$.

{\em Stage 2.} The previous stage showed that by using 
just \oone\ and \othr, the links in $M$ can be altered to the point
where they represent diffeomorphic 4-manifolds, namely the common
``middle slice'' $W_{1/2}$, henceforth renamed simply $W$. Note
that nothing in that stage perturbed the canonical decomposition of
the boundary of $W$ into $M \cup (\partial M \times I) \cup N$, since
all the handles were attached in the interior of $M$. Now apply Cerf
theory: regard the two different handle decompositions of $W$ as being
induced by Morse functions $f_0$, $f_1$ each taking the value 0 on
$M$, $t$ on $\partial M \times \{t\}$, 1 on $N$, and having no
critical points near the boundary. These may be connected by a generic
path of functions (all satisfying the same boundary conditions) in the
usual way. There is {\em no} difference at this stage between the
original proof (where, at this stage, $W$ has two disjoint closed
boundary components $M,N$, assigned the values 0 and 1) and this
refined case (where $M, N$ have boundary and there is a connecting
collar $\partial M \times I$ with projection to $I$ as a Morse
function) since Cerf's theory of genericity of functions on manifolds
with boundary (see for example \cite{HW}) requires only that there are
never critical points near the boundary: this is ensured by the
restriction of the Morse function. Consequently, everything can be
worked as in \cite{K1}. The arc of functions defines a graphic
depicting the indices and heights of the critical points; there are
isolated moments at which a handle pair is born or dies, or when one
handle's descending manifold slides over another's critical point,
instead of descending to the ``base'' $M$, but generically all that is
occurring is isotopy of the attaching maps. (Other than at
handleslides, all descending manifolds still reach the base because the
Morse function on the ``sides'' is arranged to stop them landing
there.)

Recall that the idea of Kirby's proof is to alter the graphic so as to
remove all critical points of index other than 2: for then, what
remains is a graphic depicting isotopy of 2-handle attaching maps and
occasional slides of one 2-handle over another, which correspond to
the move \otwo. Much progress can be made by applying the manouevres
of Cerf theory (dovetail lemma, beak lemma, principle of independent
trajectories). First, the 0-1-handle births and deaths are put in a
nested form. By introducing trivial 1-2-pairs and cancelling their
1-handles with the 0-handles (formally, introducing two dovetails and
then moving beaks and independent trajectories), all 0-handles can be
removed from the graphic, leaving new 2-handles in their place. A
similar procedure replaces the 4-handles by 2-handles. Next, the 1-2-
and 2-3-pairs are nested in a similar fashion: it is advantageous at
this stage to apply the above procedure once again, in order to
replace the 3-handles by 1-handles (this was not used in
\cite{K1}). So at the end of this stage, the graphic has only 1- and
2-handles.

{\em Stage 3.} There remains only the problem of ``continuously
surgering a 1-handle to a 2-handle'' all the way along the arc of
functions between the birth and death.

Consider the outermost 1-2-pair in the nest. Just after the birth, the
attaching maps may be visualised as a pair of 3-balls in $M$,
representing the feet of the 1-handle, together with a (framed)
attaching curve $A \cup C$ for the 2-handle, where $C$ is an arc in
$M$ and $A$ is the core of the 1-handle.  The surgery will be
specified by picking any (framed) curve consisting of another core
$A'$ of the 1-handle and an arc $B_0$ in $M$ (choose it to miss all
the other handle attachments). During the period of time between the
birth and death of the 1-handle under consideration, the feet of the
1-handle are isotoped around, as are the other attaching maps, and
there are births and deaths of pairs and 2-handle slides. Notice that
1-handles do {\em not} slide on 1-handles: this is part of the outcome
of the Cerf theory in stage 2. The isotopy of the feet extends to an
ambient isotopy of $M$ which drags the arc $B_0$ through a family
$B_t$.

``Continuous surgery'' on $A' \cup B_t$ means replacing its
neighbourhood by $B^2 \times S^2$ in some standard fashion,
independent of $t$.  A concrete way of realising this 5-dimensional
cobordism is to attach a pair of 2-handles (as described more
explicitly in stage 1) onto $W$, one running along $A' \cup B_t$ and
one (dual to it) on a small 0-framed meridian of this curve. Then
perform a collapse of the 1-handle across the new `long' 2-handle onto
a regular neighbourhood of $B_t$, removing the pair. The effect on the
handle decomposition is to delete the 1-handle, connect up all the
other attaching curves which run over it using arcs running parallel
to $B_t$ in its regular neighbourhood, and encircle all these with a
0-framed unknot. This operation is obviously continuous with respect
to isotopy of attaching maps (and to other births and deaths and
handleslides, which may be chosen to happen away from the feet and
arc) except when a 2-handle attaching curve crosses the arc $B_t$, but
this requires only a slide move \otwo. This procedure may thus be used
to eliminate the 1-handles, but it remains to examine the change in
framed links which occurs in the time interval spanning the birth and
the surgery. (The case of a death is exactly the reverse of this
situation, once any attaching curves which go over the 1-handle have
been slid off over the cancelling 2-handle, using move \otwo).

\begin{figure}
\[\vcenter{\epsffile{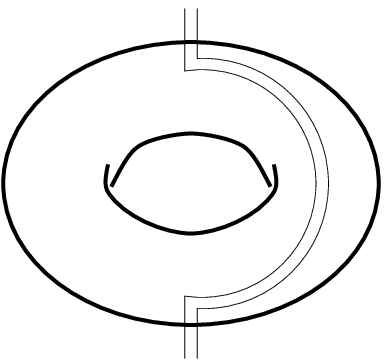}}\ \rightarrow \vcenter{\epsffile{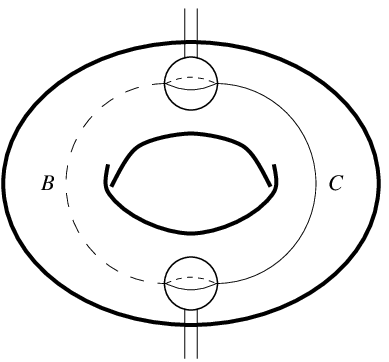}} 
\rightarrow \vcenter{\epsffile{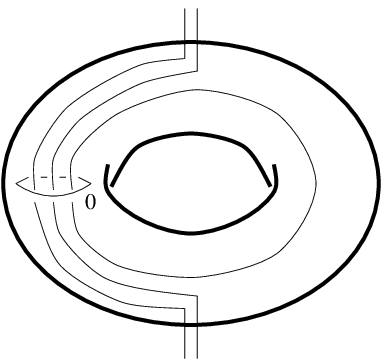}} \]
\caption{Birth of a 1-2-handle pair via expansion from $C$, followed by the surgery defined by $B$ and a collapse to its regular neighbourhood. \label{f:1-h}} 
\end{figure}

The birth of such a pair can be explicitly realised (via a
homeomorphism) as an expansion (inverse collapsing) move from a
regular neighbourhood of the arc $C \subseteq M$ (which is assumed to
be disjoint from any other attaching curves.)  The changes in the link
are shown in figure \ref{f:1-h}; the initial and final configurations
differ by handleslides and then a circumcision \othr, finishing the
proof.
\end{proof}

\begin{figure}
\[ \vcenter{\epsffile{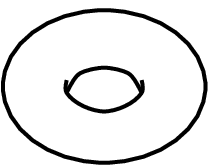}}\ \longrightarrow\ \vcenter{\epsffile{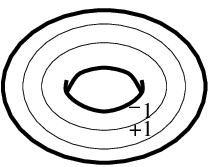}}\]
\caption{{\em Antigeny} ($\othrd$) move.\label{f:3'}}
\end{figure}

\begin{remark} It is easily shown that the moves \oone, \otwo, \othr\
generate the same equivalence relation on framed links as do \oone,
\otwo, $\othrd$, where $\othrd$ is the alternative
operation shown in figure \ref{f:3'}. 
\end{remark}

\begin{remark} In this formulation of the theorem, the question of 
``reparametrisation of the boundary'' has been deliberately
avoided. It can be brought in by gluing on a mapping cylinder, defined
as $(F \times [0, \frac12]) \cup_f (F \times [\frac12, 1])$, where $F$
is a closed surface, and $f$ is an automorphism: in this form, its
boundaries are still canonically identified. Expressing $f$ as a
product of Dehn twists gives a presentation of the cylinder as $F
\times I$, surgered on a sequence of curves in successive $F$-slices
of the cylinder. The move $\othrd$, in this context, expresses
cancellation of a Dehn twist and its inverse in the mapping class
group. (It seems to me possible, though unlikely, that Kirby calculus
in $F \times I$ could be used to derive a presentation of the mapping
class group of $F$.)
\end{remark}

\begin{center}{\em The non-orientable case}\end{center}

If $M$ is non-orientable then it is still possible to give such a
classification theorem. Once again, if $M$ and $N$ are two compact,
connected, non-orientable manifolds with identified boundaries, then
gluing them via a collar on the boundary gives a closed non-orientable
manifold. By Lickorish's theorem \cite{L2}, it bounds a (smooth)
non-orientable 4-manifold. A handle decomposition of this 4-manifold
built on $M \times I$ can have its 1- and 3-handles removed by
surgering them to 2-handles (surgery is performed on a circle
consisting of the core union an arc in $M$, which has to be chosen so
that the circle has a trivial normal bundle; since $M$ itself contains
orientation-reversing loops, this will always be possible).
Consequently, it is possible to get from $M$ to $N$ ({\em rel} the
boundary) by surgery on a framed link. It should be noted that since
2-handles attach along solid tori rather than twisted disc bundles,
the homotopy classes of all attaching curves lie in the
orientation-preserving subgroup of $\pi_1(M)$. (Note also that adding
2-handles cannot change orientability of the 3-manifold, so this is a
totally separate case.)

\begin{figure}
\[ \vcenter{\epsffile{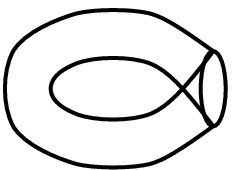}}\ \longrightarrow\ \vcenter{\epsffile{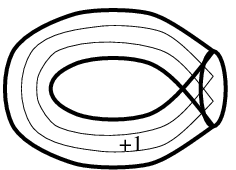}}\]
\caption{{\em M\"obius} (\ofou) move. \label{f:4}}
\end{figure}

To generate the equivalence relation on links requires another move
\ofou, similar to the $\mu$-move which was introduced by Fenn and
Rourke \cite{FR}. It is supported in a solid Klein bottle $S^1
\tilde{\times} B^2$ contained in $M$, and shown schematically in figure 
\ref{f:4}. Take the simple closed curve
on the bounding Klein bottle which runs twice around the $S^1$
direction (this is unique up to isotopy, in fact Lickorish \cite{L2}
demonstrates that there are only two non-trivial isotopy classes of
unoriented, orientation-preserving simple closed curves on a Klein
bottle).  Give this curve a framing $+1$ relative to the surface of
the Klein bottle, and push it slightly into the solid bottle. The move
\ofou\ consists of doing surgery on this framed curve. It is easy to
see that the surgered manifold is homeomorphic ({\em rel} its
boundary) to the same solid Klein bottle by considering the
equivalence between surgery on a curve in a surface (with relative
framing $+1$) and cutting-and-regluing via a (negative) Dehn twist on
that curve. Since the curve on the Klein bottle bounds a M\"obius
strip, and the Dehn twist parallel to the boundary of such a strip is
isotopic to the identity, the homeomorphism is clear. (This also makes
it clear that parity of the framing is all that matters.)

An alternative interpretation of this move at the 4-manifold level is
useful. Consider the standard decomposition of ${\mathbb R}P^4$ with one
handle in each dimension. Its 2-skeleton is $S^1 \tilde{\times} B^3$
union a 2-handle along the framed curve above. (The framing is checked
by computing the mod 2 self-intersection of ${\mathbb R}P^2$, which is
the core of the 2-handle union the M\"obius strip used above.) The
other two handles form $S^1 \tilde{\times} B^3$. An $S^1 \tilde{\times} B^3$
contained in a 4-manifold may be cut out and replaced with the
complement of the 1-skeleton of ${\mathbb R}P^4$ by using the 5-cobordism
$(S^1 \tilde{\times} B^3) \times I$. If the original $S^1
\tilde{\times} B^3$ is a 4-dimensional neighbourhood of an
orientation-reversing curve in $M \subseteq M \times I$, then
performing this replacement corresponds to attaching the 2-handle to
the curve in the 3-dimensional solid Klein bottle neighbourhood in
$M$, as in move \ofou.

\begin{theorem} If $M$ is a compact, connected,
non-orientable manifold then moves \oone, \otwo, \othr, \ofou\
suffice.
\end{theorem}

\begin{proof}
The start of the proof is the same as before: connect the 4-manifolds
via a collar. The resulting closed 4-manifold bounds a smooth
5-manifold if and only if its Stiefel-Whitney numbers $w_2^2$ and
$w_4$ to vanish. This can be achieved by connect-summing with ${\mathbb
C}P^2$ ($w_2^2=1, w_4=1$) and ${\mathbb R}P^4$ ($w_2^2=0, w_4=1$) if
necessary. The first is effected by move \oone. In the second case,
use the 5-cobordism described above instead of a 1-handle: the effect
on the characteristic classes is the same, and this is move
\ofou.

Now take a Morse function on the 5-manifold whose restriction to the
boundary is as before. Once again, replace the 1- and 4-handles by 3-
and 2-handles using surgery, having first chosen suitable arcs to get
trivial normal bundles. The effects of the 2- and 3-handles are, as
before, to add pairs of 2-handles corresponding to move \othr,
starting from each end of the 5-manifold, to reach the common
$W_{1/2}$. Now the Cerf theory works as before, and lack of
orientability does not play a part in the remainder of the proof.
\end{proof}

\noindent{\em Acknowledgement.}
Recently Kerler \cite{Ke1} and Sawin \cite{S}, motivated by the desire
to give clean constructions of Topological Quantum Field Theories \`a
la Reshetikhin-Turaev, obtained related presentations and `moves' for
general orientable 3-manifolds. Neither of these papers dealt with the
natural generalisation of the original framed link calculus of Kirby,
as presented here, although Kerler checked in \cite{Ke2} that his
bridged link calculus (essentially, allowing 1-handles as well as
2-handles) implied it. (The rather different version of Sawin does not
seem to have this property.) Neither did these papers deal with the
non-orientable case (usually neglected in TQFT), although presumably
they would generalise too.  I am grateful to Rob Kirby and Steve Sawin
for the encouragement to write this version up.


\begin{thebibliography}{Ke2}
\bibitem[FR]{FR} R. A. Fenn, C. P. Rourke. {\em On Kirby's calculus of
links}, Topology 18 (1979), 1-15.
\bibitem[HW]{HW} A. E. Hatcher, J. B. Wagoner. {\em Pseudo-isotopies of
 compact manifolds}, Ast\'erisque 6 (1973).
\bibitem[Ke1]{Ke1} T. Kerler. {\em Bridged links and tangle
presentations of cobordism categories}, Adv.\ Math. (to appear).
\bibitem[Ke2]{Ke2} T. Kerler. {\em Equivalence of a bridged link calculus
and Kirby's calculus of links on non-simply connected 3-manifolds},
preprint (1994).
\bibitem[K1]{K1} R. C. Kirby. {\em A calculus for framed links in $S^3$},
 Invent.\ Math.\ 45 (1978), 35-56.
\bibitem[K2]{K2} R. C. Kirby. {\em Topology of 4-manifolds}, Springer
Lect. Notes in Maths. 1374 (1989).
\bibitem[L1]{L1} W. B. R. Lickorish. {\em A representation of orientable
combinatorial 3-manifolds}, Ann.\ Math.\ 76 (1962), 531-540.
\bibitem[L2]{L2} W. B. R. Lickorish. {\em Homeomorphisms of
non-orientable 2-manifolds}, Proc.\ Cam.\ Phil.\ Soc.\ 59 (1963), 307-317.
\bibitem[S]{S} S. Sawin. {\em Extending the Kirby calculus to
manifolds with boundary}, preprint (1994).
\end{thebibliography}
\end{document}